\documentclass[12pt]{article}
\usepackage{amsfonts}
\usepackage{amsmath}
\usepackage{graphicx}

\input tcilatex
\input tcilatex
\begin{document}

\author{Ozlem Ersoy\thanks{%
Ozlem Ersoy, Telephone: +90 222 2393750, Fax: +90 222 2393578, E-mail:
ozersoy@ogu.edu.tr}, Idris Dag and Nihat Adar \\
Eski\c{s}ehir Osmangazi University, \\
Faculty of Science and Art,\\
Department of Mathematics-Computer, Eski\c{s}ehir, Turkey,}
\title{\textbf{The Exponential Cubic B-spline Algorithm for Burgers' Equation%
}}
\maketitle

\begin{abstract}
\textbf{Purpose --} The exponential cubic B-spline functions are used to set
up the collocation method for finding solutions of the Burgers' equation.
The effect of the exponential cubic B-splines in the collocation method is
sought by studying four text problems.

\textbf{Design --} The Burgers' equation is fully-discretized using the
Crank-Nicolson method for the time discretizion and exponential cubic
B-spline function for discretizion of spatial variable.

\textbf{Findings --} The exponential cubic B-spline methods have produced
acceptable solution of Burgers' equation.

\textbf{Originality/value --} Burgers' equation has never been solved by the
collocation method using exponential cubic B-spline.


\end{abstract}

\section{Introduction}

This paper is concerned with adapting the exponential cubic B-spline
function into the collocation method to develop a numerical method for
finding numerical solutions of the Burgers' equation of the form%
\begin{equation}
U_{t}+UU_{x}-\lambda U_{xx}=0,\text{ }a\leq x\leq b,\text{ }t\geq 0
\label{1}
\end{equation}%
with the initial condition%
\begin{equation}
U(x,0)=f(x),\text{ }a\leq x\leq b  \label{bas}
\end{equation}%
and the boundary conditions%
\begin{equation}
U(a,t)=\sigma _{1},U(b,t)=\sigma _{2}  \label{sin}
\end{equation}%
where subscripts $x$ and $t$ denote differentiation, $\lambda =\frac{1}{%
\func{Re}}>0$ and $\func{Re}$ is the Reynolds number characterizing the
strength of viscosity, $\sigma _{1},\sigma _{2}$ are the constants $u=u(x,t)$
is a sufficiently differentiable unknown function and $f(x)$ is a bounded
function. Initial condition and boundary conditions will be defined in the
later section depending on the test problems.

Burgers' equation (\ref{1}) has been widely used in modelling various
problem in science and engineering with a broad range of applications,
including gas dynamics, traffic flow, wave propagation in acoustics and
hydrodynamics, etc. It is solved exactly for an arbitrary initial condition 
\cite{ho} so that it has been used as an test equation for the numerical
methods. But some exact solutions contain series solution and these
solutions are unpractical because of the slow convergence of the Fourier
series when the small values of viscosity constant are used. Thus the
numerical studies are necessary to find improsing solutions of the Burgers'
equation. Besides the Burgers' equation with small viscosity constant gives
rises to appearance of steep front and shock waves in the solutions.
Therefore one comes across complexity in writing the influential numerical
method for finding solutions of the Burgers' equation. Many numerical
studies have been published on the numerical solution of the Burgers'
equation. Especially, attention has been given in deriving the numerical
methods for the numerical solution of the Burgers' equation with small
values of viscosity constant.

Burgers' equation was first introduced by \cite{bat}. Solutions of the
Burgers' equation were presented by using some numerical methods with
splines. A cubic spline collocation procedure has been developed for the
numerical solution of the Burgers' equation in the papers\cite{rubin, ru1} .
A B-spline Galerkin Method is described to solve the Burgers' equation over
the both fixed and varied distribution of knots to define the B-splines in
the studies of Davies \cite{ad,ad1}.\ A Numerical method \cite{u1,nu2, nu3}
is developed \ for solving the Burgers' equation by using splitting method
and cubic spline approximation method. In \cite%
{er,ras,bs,irk,ram,bbs,ofem2,bss1,rr,cv}, numerical solutions of the
one-dimensional Burgers' equation are obtained by a methods based on
collocation of quadratic, cubic, quintic and septic B-splines over finite
elements respectively, in which approximate functions in the collocation
method for the Burgers' equation are constructed by using the various degree
B-splines. Galerkin methods based on various degree B-splines have been set
up to find approximate solutions of the Burgers' equation in the studies 
\cite{ad1,ofem3, bbs}. The least square method is combined with the
B-splines to form numerical methods for solving the Burgers' equation in the
works\cite{sku,ww}. Numerical solutions of the Burgers' equation is
presented based on the cubic B-spline quasi-interpolation and the compact
finite difference method in \cite{zi,re}. Taylor-collocation and
Taylor-Galerkin methods for the numerical solutions of the Burgers' equation
are formed by using both cubic and quadratic B-splines in the study \cite%
{ofem8} respectively. \ Differential quadrature methods based on cubic and
quartic B-splines are set up to solve the Burgers equation in the works \cite%
{al1,al,yu}. The hybrid spline difference method is developed to solve the
Burgers' equation by Chi-Chang Wang and his colleagues\cite{hi}. The
nonlinear polynomial splines method been proposed for solving the Burgers
equation by M. Caglar and M. F. Ucar recently.

The exponential cubic B-spline function and its some properties are
described in detail in the paper\cite{maccartin}. Since each exponential
basis we use is twice continuously differentiable, we can form twice
continuously approximate solution \ to the differential equations. There
exist few articles which are used to form numerical methods to solve
differential equations. The exponential cubic B-splines are used with the
collocation method to find the numerical solution of the singular
perturbation problem by Manabu \ Sakai and his colleague\cite{sc3}. Another
application of the collocation method using the cardinal exponential cubic
B-splines was shown for finding the numerical solutions of the singularly
perturbed boundary problem in the study of Desanka Radunuvic \cite{sc4}. The
exponential cubic B-spline collocation method is set up to obtain the
numerical solutions of the self-adjoint singularly perturbed boundary value
problems in the work\cite{sch}. The only linear partial differential
equation known as the convection-diffusion equation is solved by way of the
exponential cubic B-spline collocation method in the study\cite{sc2}.

In this paper, we have compared results \ of the Burgers' equation with
those obtained with both the cubic B-spline collocation method and cubic
B-spline Galerkin finite element method \cite{bs,bbs} since the B-spline and
exponential cubic B-spline functions have almost the same properties. In
section 2, exponential cubic B-spline collocation method is described. In
section 3, four classical test examples are studied to show the versatility
of proposed algorithm and finally the conclusion is included to discuss \
the outcomes of the algorithm.

\section{Exponential Cubic B-spline Collocation Method}

Knots are equally distributed over the problem domain as 
\begin{equation*}
\pi :a=x_{0}<x_{1}<\ldots <x_{N}=b
\end{equation*}%
with mesh spacing $h=(b-a)/N.$ The exponential cubic B-splines, $B_{i}${}$%
(x),$ at the points of $\pi ,$ can be defined as

\begin{equation}
B_{i}(x)=\left \{ 
\begin{array}{ll}
b_{2}\left( \left( x_{i-2}-x\right) -\dfrac{1}{p}\left( \sinh (p\left(
x_{i-2}-x\right) )\right) \right) & \left[ x_{i-2},x_{i-1}\right] , \\ 
a_{1}+b_{1}\left( x_{i}-x\right) +c_{1}\exp \left( p\left( x_{i}-x\right)
\right) +d_{1}\exp \left( -p\left( x_{i}-x\right) \right) & \left[
x_{i-1},x_{i}\right] , \\ 
a_{1}+b_{1}(x-x_{i})+c_{1}\exp \left( p\left( x-x_{i}\right) \right)
+d_{1}\exp \left( -p\left( x-x_{i}\right) \right) & \left[ x_{i},x_{i+1}%
\right] , \\ 
b_{2}\left( (x-x_{i+2})-\dfrac{1}{p}(\sinh \left( p\left( x-x_{i+2}\right)
\right) )\right) & \left[ x_{i+1},x_{i+2}\right] , \\ 
0 & \text{otherwise.}%
\end{array}%
\right.  \label{r3}
\end{equation}%
where%
\begin{equation*}
\begin{array}{l}
a_{1}=\dfrac{phc}{phc-s},\text{ }b_{1}=\dfrac{p}{2}[\dfrac{c(c-1)+s^{2}}{%
(phc-s)(1-c)}],\text{ }b_{2}=\dfrac{p}{2(phc-s)}, \\ 
c_{1}=\dfrac{1}{4}[\dfrac{\exp (-ph)(1-c)+s(\exp (-ph)-1)}{(phc-s)(1-c)}] \\ 
d_{1}=\dfrac{1}{4}[\dfrac{\exp (ph)(c-1)+s(\exp (ph)-1)}{(phc-s)(1-c)}]%
\end{array}%
\end{equation*}%
and $c=\cosh (ph),$ $s=\sinh (ph),p$ is a free parameter. When $p=1$, graph
of the exponential cubic B-splines over the interval $[0.1]$ is depicted in
Fig. 1.%
\begin{equation*}
\begin{array}{c}
\FRAME{itbpF}{2.4431in}{2.4431in}{0in}{}{}{fig1.bmp}{\special{language
"Scientific Word";type "GRAPHIC";maintain-aspect-ratio TRUE;display
"USEDEF";valid_file "F";width 2.4431in;height 2.4431in;depth
0in;original-width 3.2197in;original-height 3.2197in;cropleft "0";croptop
"1";cropright "1";cropbottom "0";filename 'fig/Fig1.bmp';file-properties
"XNPEU";}} \\ 
\text{Fig.1: exponential cubic B-splines over the interval }[0.1]%
\end{array}%
\end{equation*}

$\{B_{-1}(x),B_{0}(x),\cdots ,B_{N+1}(x)\}$ forms a basis for the functions
defined over the interval. Each basis function $B_{i}(x)$ is twice
continuously differentiable. The values of $B_{i}(x),B_{i}^{^{\prime }}(x)$
and $B_{i}^{^{\prime \prime }}(x)$ at the knots $x_{i}$ 's can be computed
from Eq.(\ref{r3}) shown Table 1.

\begin{equation*}
\begin{tabular}{l}
Table 1: Values of $B_{i}(x)$ and its principle two \\ 
derivatives at the knot points \\ 
$%
\begin{tabular}{|l|ccccc|}
\hline
$x$ & $x_{i-2}$ & $x_{i-1}$ & $x_{i}$ & $x_{i+1}$ & $x_{i+2}$ \\ \hline
$B_{i}$ & $0$ & $\dfrac{s-ph}{2(phc-s)}$ & $1$ & $\dfrac{s-ph}{2(phc-s)}$ & $%
0$ \\ 
$B_{i}^{^{\prime }}$ & $0$ & $\dfrac{p(1-c)}{2(phc-s)}$ & $0$ & $\dfrac{%
p(c-1)}{2(phc-s)}$ & $0$ \\ 
$B_{i}^{^{\prime \prime }}$ & $0$ & $\dfrac{p^{2}s}{2(phc-s)}$ & $-\dfrac{%
p^{2}s}{phc-s}$ & $\dfrac{p^{2}s}{2(phc-s)}$ & $0$ \\ \hline
\end{tabular}%
$%
\end{tabular}%
\end{equation*}

Now suppose that an approximate solution $U_{N}$ to the unknown $U$ is given
by

\begin{equation}
U_{N}(x,t)=\sum_{i=-1}^{N+1}\delta _{i}B_{i}(x)  \label{r4}
\end{equation}%
where $\delta _{i}$ are time dependent parameters to be determined from the
collocation points $x_{i},i=0,...,N$, the boundary and initial conditions.
Evaluation of Eq. (\ref{r4}), its first and second derivatives at knots $%
x_{i}$ using the Table 1 yields the nodal values U$_{i}$ in terms of
parameters%
\begin{equation}
\begin{tabular}{l}
$U_{i}=U(x_{i},t)=\dfrac{s-ph}{2(phc-s)}\delta _{i-1}+\delta _{i}+\dfrac{s-ph%
}{2(phc-s)}\delta _{i+1},$ \\ 
$U_{i}^{\prime }=U^{\prime }(x_{i},t)=\dfrac{p(1-c)}{2(phc-s)}\delta _{i-1}+%
\dfrac{p(c-1)}{2(phc-s)}\delta _{i+1}$ \\ 
$U_{i}^{\prime \prime }=U^{\prime \prime }(x_{i},t)=\dfrac{p^{2}s}{2(phc-s)}%
\delta _{i-1}-\dfrac{p^{2}s}{phc-s}\delta _{i}+\dfrac{p^{2}s}{2(phc-s)}%
\delta _{i+1}.$%
\end{tabular}
\label{r5}
\end{equation}

The Crank--Nicolson \ scheme is used to discretize time variables of the
unknown U in the Burgers' equation so that one obtain the time discretized
form of the equation as 
\begin{equation}
\frac{U^{n+1}-U^{n}}{\Delta t}+\frac{(UU_{x})^{n+1}+(UU_{x})^{n}}{2}-\lambda 
\frac{U_{xx}^{n+1}+U_{xx}^{n}}{2}=0  \label{r6}
\end{equation}%
where $U^{n+1}=U(x,t)$ is the solution of the equation at the $(n+1)$th time
level. Here $t^{n+1}$ $=t^{n}+t$, and $\Delta t$ is the time step,
superscripts denote $n$ th time level , $t^{n}=n\Delta t$

The nonlinear term $(UU_{x})^{n+1}$ in Eq. (\ref{r6}) is linearized by using
the following form \cite{rubin}:%
\begin{equation}
(UUx)^{n+1}=U^{n+1}U_{x}^{n}+U^{n}U_{x}^{n+1}-U^{n}U_{x}^{n}  \label{r7}
\end{equation}%
So Equation (\ref{r6}) is discretized in time as%
\begin{equation}
U^{n+1}-U^{n}+\frac{\Delta t}{2}(U^{n+1}U_{x}^{n}+U^{n}U_{x}^{n+1})-\lambda 
\frac{\Delta t}{2}(U_{xx}^{n+1}-U_{xx}^{n})=0  \label{r8}
\end{equation}%
Substitution of \ref{r4} into \ref{r8} leads to the fully-discretized
equation:

\begin{equation}
\begin{tabular}{l}
$\left( \alpha _{1}+\dfrac{\Delta t}{2}\left( \alpha _{1}L_{2}+\beta
_{1}L_{1}-\lambda \gamma _{1}\right) \right) \delta _{m-1}^{n+1}+\left(
\alpha _{2}+\dfrac{\Delta t}{2}\left( \alpha _{2}L_{2}-\lambda \gamma
_{2}\right) \right) \delta _{m}^{n+1}+$ \\ 
$\left( \alpha _{3}+\dfrac{\Delta t}{2}\left( \alpha _{3}L_{2}+\beta
_{2}L_{1}-\lambda \gamma _{3}\right) \right) \delta _{m+1}^{n+1}=(\alpha
_{1}+\lambda \dfrac{\Delta t}{2}\gamma _{1})\delta _{m-1}^{n}+$ \\ 
$(\alpha _{2}+\lambda \dfrac{\Delta t}{2}\gamma _{2})\delta _{m}^{n}+(\alpha
_{3}+\lambda \dfrac{\Delta t}{2}\gamma _{3})\delta _{m+1}^{n}$%
\end{tabular}
\label{r9}
\end{equation}%
where%
\begin{eqnarray*}
L_{1} &=&\alpha _{1}\delta _{i-1}+\alpha _{2}\delta _{i}+\alpha _{3}\delta
_{i+1} \\
L_{2} &=&\beta _{1}\delta _{i-1}+\beta _{2}\delta _{i+1}
\end{eqnarray*}%
\begin{eqnarray*}
\alpha _{1} &=&\dfrac{s-ph}{2(phc-s)},\text{ }\alpha _{2}=1,\text{ }\alpha
_{3}=\dfrac{s-ph}{2(phc-s)} \\
\beta _{1} &=&\dfrac{p(1-c)}{2(phc-s)},\text{ }\beta _{2}=\dfrac{p(c-1)}{%
2(phc-s)} \\
\gamma _{1} &=&\dfrac{p^{2}s}{2(phc-s)},\text{ }\gamma _{2}=-\dfrac{p^{2}s}{%
phc-s},\text{ }\gamma _{3}=\dfrac{p^{2}s}{2(phc-s)}
\end{eqnarray*}

The system consist of $N+1$ linear equation in $N+3$ unknown parameters $%
\mathbf{d}^{n+1}=(\delta _{-1}^{n+1},\delta _{0}^{n+1},\ldots ,\delta
_{N+1}^{n+1})$. To make solvable the system, boundary conditions $\sigma
_{1}=U_{0},$ $\sigma _{2}=U_{N}$ are used to find two additional linear
equations: 
\begin{eqnarray}
\delta _{-1} &=&\frac{1}{\alpha _{1}}\left( U_{0}-\alpha _{2}\delta
_{0}-\alpha _{3}\delta _{1}\right) ,  \label{r10} \\
\delta _{N+1} &=&\frac{1}{\alpha _{3}}\left( U_{N}-\alpha _{1}\delta
_{N-1}-\alpha _{2}\delta _{N}\right) .  \notag
\end{eqnarray}%
(\ref{r10}) can be used to eliminate $\delta _{-1},\delta _{N+1}$ from the
system (\ref{r9}) which then becomes the solvable matrix equation for the
unknown $\delta _{0}^{n+1},\ldots ,\delta _{N}^{n+1}.$ A variant of Thomas
algorithm is used to solve the system.

Initial parameters $\delta _{-1}^{0},\delta _{0}^{0},\ldots ,\delta
_{N+1}^{0}$ can be determined from the initial condition and first space
derivative of the initial conditions at the boundaries as the following:

\begin{enumerate}
\item $U_{N}(x_{i},0)$ $=U(x_{i},0),$ $i=0,...,N$

\item $(U_{x})_{N}(x_{0},0)=U^{\prime }(x_{0})$

\item $(U_{x})_{N}(x_{N},0)=U^{\prime }(x_{N}).$
\end{enumerate}

\section{Numerical tests}

\qquad Numerical method described in the previous section will be tested on
three text problems for getting solutions of \ the Burgers' equation. Four
kinds of examples are presented in order to demonstrate the versatility\ and
the accuracy of the proposed method.\ The discrete $L_{2}$ and $L_{\infty }$
error norm%
\begin{equation*}
\begin{array}{l}
L_{2}=\left \vert U-U_{N}\right \vert =\sqrt{h\sum \limits_{j=0}^{N}\left
\vert (U_{j}-(U_{N})_{j}^{n})^{2}\right \vert }, \\ 
L_{\infty }=\left \vert U-U_{N}\right \vert _{\infty }=\max \limits_{j}\left
\vert U_{j}-(U_{N})_{j}^{n}\right \vert%
\end{array}%
\end{equation*}%
are used to measure error between the analytical and numerical solutions.

\textbf{(a)} The Burger's equation, with the sine wave initial condition $%
U(x,0)=sin(\pi x)$ and boundary conditions $U(0,t)=U(1,t)=0,$ has analytic
solution in the form of the infinite series defined by \cite{cole}\ as

\begin{equation}
U(x,t)=\frac{4\pi \lambda \dsum \limits_{j=1}^{\infty }j\mathbf{I}_{j}(\frac{%
1}{2\pi \lambda })\sin (j\pi x)\exp (-j^{2}\pi ^{2}\lambda t)}{\mathbf{I}%
_{0}(\frac{1}{2\pi \lambda })+2\dsum \limits_{j=1}^{\infty }\mathbf{I}_{j}(%
\frac{1}{2\pi \lambda })\cos (j\pi x)\exp (-j^{2}\pi ^{2}\lambda t)}
\end{equation}%
where $\mathbf{I}_{j}$ are the modified Bessel functions. This problem gives
the decay of sinusoidal disturbance. Numerical solutions at different times
are depicted in the Figs 2-5 for the parameters $N=40,$ $\Delta
t=0.0001,\lambda =1,0.1,0.01,0.001$ from the figures we see that the smaller
viscosities $\  \lambda $ cause to develop the sharp front thorough the right
boundary and amplitude of the sharp front starts to decay as time progress.
These properties of \ solutions are in very good agreement with findings of
B Saka and \ I Dag \cite{bss, bss1}.%
\begin{equation*}
\begin{array}{cc}
\begin{tabular}{l}
\FRAME{itbpF}{2.8738in}{2.4275in}{0in}{}{}{fig2.bmp}{\special{language
"Scientific Word";type "GRAPHIC";display "USEDEF";valid_file "F";width
2.8738in;height 2.4275in;depth 0in;original-width 14.2288in;original-height
8.0004in;cropleft "0";croptop "1";cropright "1";cropbottom "0";filename
'Fig2.bmp';file-properties "XNPEU";}} \\ 
Fig. 2: Solutions at different times \\ 
for $\lambda =1,$ $N=40,$ $\Delta t=0.0001.$%
\end{tabular}
& 
\begin{tabular}{l}
\FRAME{itbpF}{2.8729in}{2.4275in}{0in}{}{}{fig3.bmp}{\special{language
"Scientific Word";type "GRAPHIC";display "USEDEF";valid_file "F";width
2.8729in;height 2.4275in;depth 0in;original-width 14.2288in;original-height
8.0004in;cropleft "0";croptop "1";cropright "0.9998";cropbottom "0";filename
'Fig3.bmp';file-properties "XNPEU";}} \\ 
Fig. 3: Solutions at different times \\ 
for $\lambda =0.1,$ $N=40,$ $\Delta t=0.0001.$%
\end{tabular}%
\end{array}%
\end{equation*}%
\begin{equation*}
\begin{array}{cc}
\begin{tabular}{l}
\FRAME{itbpF}{2.8738in}{2.4275in}{0in}{}{}{fig4.bmp}{\special{language
"Scientific Word";type "GRAPHIC";display "USEDEF";valid_file "F";width
2.8738in;height 2.4275in;depth 0in;original-width 14.2288in;original-height
8.0004in;cropleft "0";croptop "1";cropright "1";cropbottom "0";filename
'Fig4.bmp';file-properties "XNPEU";}} \\ 
Fig. 4: Solutions at different times \\ 
for $\lambda =0.01,$ $N=40,$ $\Delta t=0.0001.$%
\end{tabular}
& 
\begin{tabular}{l}
\FRAME{itbpF}{2.8738in}{2.4275in}{0in}{}{}{fig5.bmp}{\special{language
"Scientific Word";type "GRAPHIC";display "USEDEF";valid_file "F";width
2.8738in;height 2.4275in;depth 0in;original-width 14.2288in;original-height
8.0004in;cropleft "0";croptop "1";cropright "1";cropbottom "0";filename
'Fig5.bmp';file-properties "XNPEU";}} \\ 
Fig. 5: Solutions at different times \\ 
for $\lambda =0.001,$ $N=40,$ $\Delta t=0.0001.$%
\end{tabular}%
\end{array}%
\end{equation*}%
Two dimensional solutions are depicted from time $t=0$ to $t=1$ with time
increment $\Delta t=0.01$ for space increment $h=0.25$ and various $\lambda $
in Fig. 6-9. When the smaller $\lambda =0.001$ is taken, the solutions
starts to decay after about time $t=0.6$ when $N=40$ is used. So to have
acceptable solution with $\lambda =0.001,$we decrease the space step to $%
h=0.02$ and graph of the solution is shown in Fig 9.

\begin{equation*}
\begin{tabular}{cc}
\FRAME{itbpF}{3.0193in}{2.8206in}{0in}{}{}{fig6.bmp}{\special{language
"Scientific Word";type "GRAPHIC";maintain-aspect-ratio TRUE;display
"USEDEF";valid_file "F";width 3.0193in;height 2.8206in;depth
0in;original-width 5.156in;original-height 4.8136in;cropleft "0";croptop
"1";cropright "1";cropbottom "0";filename 'fig/Fig6.bmp';file-properties
"XNPEU";}} & \FRAME{itbpF}{3.1514in}{2.8426in}{0in}{}{}{fig7.bmp}{\special%
{language "Scientific Word";type "GRAPHIC";maintain-aspect-ratio
TRUE;display "USEDEF";valid_file "F";width 3.1514in;height 2.8426in;depth
0in;original-width 5.2062in;original-height 4.6933in;cropleft "0";croptop
"1";cropright "1";cropbottom "0";filename 'fig/Fig7.bmp';file-properties
"XNPEU";}} \\ 
Fig. 6: Solutions for $\lambda =1,$ $N=40,$ $\Delta t=0.01$ & Fig. 7:
Solutions for $\lambda =0.1,$ $N=40,$ $\Delta t=0.01$ \\ 
\FRAME{itbpF}{3.1514in}{2.8349in}{0in}{}{}{fig8.bmp}{\special{language
"Scientific Word";type "GRAPHIC";maintain-aspect-ratio TRUE;display
"USEDEF";valid_file "F";width 3.1514in;height 2.8349in;depth
0in;original-width 5.2062in;original-height 4.6933in;cropleft "0";croptop
"1";cropright "1";cropbottom "0.0029";filename
'fig/Fig8.bmp';file-properties "XNPEU";}} & \FRAME{itbpF}{3.0193in}{2.8206in%
}{0in}{}{}{fig9.bmp}{\special{language "Scientific Word";type
"GRAPHIC";maintain-aspect-ratio TRUE;display "USEDEF";valid_file "F";width
3.0193in;height 2.8206in;depth 0in;original-width 5.156in;original-height
4.8136in;cropleft "0";croptop "1";cropright "1";cropbottom "0";filename
'fig/Fig9.bmp';file-properties "XNPEU";}} \\ 
Fig. 8: Solutions for $\lambda =0.01,$ $N=40,$ $\Delta t=0.01.$ & Fig. 9:
Solutions for $\lambda =0.001,$ $N=50,$ $\Delta t=0.001.$%
\end{tabular}%
\end{equation*}

A comparison has been made between the present collocation method \ and
alternative \ approaches including the cubic B-spline collocation method and
Cubic B-spline Galerkin finite element method for parameters of $\Delta
t=0.0001,$ $N=80,$ $\lambda =1,0.1,$ $0.01.$ Exact solutions for $\lambda
>10^{-2}$ are not practical because of the low convergence of the infinite
series so that these results are not compared with the exact solutions. It
can be seen \ from Table 2 a, 2 b, 2 c and 3 a, 3 b that accuracy of the
presented solutions is much the same with both the cubic B-spline
collocation method and cubic B-spline Galerkin finite element method. When
the size of the space variable is reduced, the error becomes less than \
that of the cubic B-spline collocation methods and is almost close to the
cubic B-spline Galerkin finite element method and solution values are
documented in Table 3\ at time $t=0.1$%
\begin{equation*}
\begin{tabular}{|c|}
\hline
$\text{Table 2 a: Comparison of the numerical solutions of Problem 1 for }%
\lambda =1$ \\ \hline
and $N=80,$ $\Delta t=0.0001$ at different times with the exact solutions \\ 
\hline
$%
\begin{tabular}{llllll}
$x$ & $t$ & $%
\begin{array}{c}
\text{Present} \\ 
p=1%
\end{array}%
$ & $%
\begin{array}{c}
\text{Ref.\cite{bs}} \\ 
\text{ }%
\end{array}%
$ & $%
\begin{array}{c}
\text{Ref \cite{bbs}} \\ 
\text{ }%
\end{array}%
$ & $%
\begin{array}{c}
\text{Exact} \\ 
\text{ }%
\end{array}%
$ \\ \hline
0.25 & 0.4 & \multicolumn{1}{c}{$0.01356$} & \multicolumn{1}{c}{$0.01357$} & 
\multicolumn{1}{c}{$0.01357$} & \multicolumn{1}{c}{$0.01357$} \\ 
& 0.6 & \multicolumn{1}{c}{$0.00189$} & \multicolumn{1}{c}{$0.00189$} & 
\multicolumn{1}{c}{$0.00189$} & \multicolumn{1}{c}{$0.00189$} \\ 
& 0.8 & \multicolumn{1}{c}{$0.00026$} & \multicolumn{1}{c}{$0.00026$} & 
\multicolumn{1}{c}{$0.00026$} & \multicolumn{1}{c}{$0.00026$} \\ 
& 1.0 & \multicolumn{1}{c}{$0.00004$} & \multicolumn{1}{c}{$0.00004$} & 
\multicolumn{1}{c}{$0.00004$} & \multicolumn{1}{c}{$0.00004$} \\ 
& 3.0 & \multicolumn{1}{c}{$0.00000$} & \multicolumn{1}{c}{$0.00000$} & 
\multicolumn{1}{c}{$0.00000$} & \multicolumn{1}{c}{$0.00000$} \\ 
&  & \multicolumn{1}{c}{} & \multicolumn{1}{c}{} & \multicolumn{1}{c}{} & 
\multicolumn{1}{c}{} \\ 
0.50 & 0.4 & \multicolumn{1}{c}{$0.01922$} & \multicolumn{1}{c}{$0.01923$} & 
\multicolumn{1}{c}{$0.01924$} & \multicolumn{1}{c}{$0.01924$} \\ 
& 0.6 & \multicolumn{1}{c}{$0.00267$} & \multicolumn{1}{c}{$0.00267$} & 
\multicolumn{1}{c}{$0.00267$} & \multicolumn{1}{c}{$0.00267$} \\ 
& 0.8 & \multicolumn{1}{c}{$0.00037$} & \multicolumn{1}{c}{$0.00037$} & 
\multicolumn{1}{c}{$0.00037$} & \multicolumn{1}{c}{$0.00037$} \\ 
& 1.0 & \multicolumn{1}{c}{$0.00005$} & \multicolumn{1}{c}{$0.00005$} & 
\multicolumn{1}{c}{$0.00005$} & \multicolumn{1}{c}{$0.00005$} \\ 
& 3.0 & \multicolumn{1}{c}{$0.00000$} & \multicolumn{1}{c}{$0.00000$} & 
\multicolumn{1}{c}{$0.00000$} & \multicolumn{1}{c}{$0.00000$} \\ 
&  & \multicolumn{1}{c}{} & \multicolumn{1}{c}{} & \multicolumn{1}{c}{} & 
\multicolumn{1}{c}{} \\ 
0.75 & 0.4 & \multicolumn{1}{c}{$0.01362$} & \multicolumn{1}{c}{$0.01362$} & 
\multicolumn{1}{c}{$0.01363$} & \multicolumn{1}{c}{$0.01363$} \\ 
& 0.6 & \multicolumn{1}{c}{$0.00189$} & \multicolumn{1}{c}{$0.00189$} & 
\multicolumn{1}{c}{$0.00189$} & \multicolumn{1}{c}{$0.00189$} \\ 
& 0.8 & \multicolumn{1}{c}{$0.00026$} & \multicolumn{1}{c}{$0.00026$} & 
\multicolumn{1}{c}{$0.00026$} & \multicolumn{1}{c}{$0.00026$} \\ 
& 1.0 & \multicolumn{1}{c}{$0.00004$} & \multicolumn{1}{c}{$0.00004$} & 
\multicolumn{1}{c}{$0.00004$} & \multicolumn{1}{c}{$0.00004$} \\ 
& 3.0 & \multicolumn{1}{c}{$0.00000$} & \multicolumn{1}{c}{$0.00000$} & 
\multicolumn{1}{c}{$0.00000$} & \multicolumn{1}{c}{$0.00000$}%
\end{tabular}%
$ \\ \hline
\end{tabular}%
\end{equation*}%
\begin{equation*}
\begin{tabular}{|c|}
\hline
$\text{Table 2 b: Comparison of the numerical solutions of Problem 1 for }%
\lambda =0.1$ \\ \hline
and $N=40,$ $\Delta t=0.0001$ at different times with the exact solutions \\ 
\hline
$%
\begin{tabular}{lllllll}
$x$ & $t$ & $%
\begin{array}{c}
\text{Present} \\ 
p=1%
\end{array}%
$ & $%
\begin{array}{c}
\text{Ref.\cite{bs}} \\ 
\text{ }%
\end{array}%
$ & $%
\begin{array}{c}
\text{Ref.\cite{ofem2}} \\ 
\text{ }%
\end{array}%
$ & $%
\begin{array}{c}
\text{Ref.\cite{bbs}} \\ 
\text{ }%
\end{array}%
$ & $%
\begin{array}{c}
\text{Exact} \\ 
\text{ }%
\end{array}%
$ \\ \hline
0.25 & 0.4 & \multicolumn{1}{c}{$0.30890$} & \multicolumn{1}{c}{$0.30890$} & 
\multicolumn{1}{c}{$0.30891$} & \multicolumn{1}{c}{$0.30890$} & 
\multicolumn{1}{c}{$0.30889$} \\ 
& 0.6 & \multicolumn{1}{c}{$0.24075$} & \multicolumn{1}{c}{$0.24075$} & 
\multicolumn{1}{c}{$0.24075$} & \multicolumn{1}{c}{$0.24074$} & 
\multicolumn{1}{c}{$0.24074$} \\ 
& 0.8 & \multicolumn{1}{c}{$0.19569$} & \multicolumn{1}{c}{$0.19569$} & 
\multicolumn{1}{c}{$0.19568$} & \multicolumn{1}{c}{$0.19568$} & 
\multicolumn{1}{c}{$0.19568$} \\ 
& 1.0 & \multicolumn{1}{c}{$0.16257$} & \multicolumn{1}{c}{$0.16258$} & 
\multicolumn{1}{c}{$0.16257$} & \multicolumn{1}{c}{$0.16257$} & 
\multicolumn{1}{c}{$0.16256$} \\ 
& 3.0 & \multicolumn{1}{c}{$0.02720$} & \multicolumn{1}{c}{$0.02720$} & 
\multicolumn{1}{c}{$0.02721$} & \multicolumn{1}{c}{$0.02720$} & 
\multicolumn{1}{c}{$0.02720$} \\ 
&  & \multicolumn{1}{c}{} & \multicolumn{1}{c}{} & \multicolumn{1}{c}{} & 
\multicolumn{1}{c}{} & \multicolumn{1}{c}{} \\ 
0.50 & 0.4 & \multicolumn{1}{c}{$0.56965$} & \multicolumn{1}{c}{$0.56965$} & 
\multicolumn{1}{c}{$0.56969$} & \multicolumn{1}{c}{$0.56964$} & 
\multicolumn{1}{c}{$0.56963$} \\ 
& 0.6 & \multicolumn{1}{c}{$0.44722$} & \multicolumn{1}{c}{$0.44723$} & 
\multicolumn{1}{c}{$0.44723$} & \multicolumn{1}{c}{$0.44721$} & 
\multicolumn{1}{c}{$0.44721$} \\ 
& 0.8 & \multicolumn{1}{c}{$0.35925$} & \multicolumn{1}{c}{$0.35925$} & 
\multicolumn{1}{c}{$0.35926$} & \multicolumn{1}{c}{$0.35924$} & 
\multicolumn{1}{c}{$0.35924$} \\ 
& 1.0 & \multicolumn{1}{c}{$0.29192$} & \multicolumn{1}{c}{$0.29192$} & 
\multicolumn{1}{c}{$0.29193$} & \multicolumn{1}{c}{$0.29191$} & 
\multicolumn{1}{c}{$0.29192$} \\ 
& 3.0 & \multicolumn{1}{c}{$0.04019$} & \multicolumn{1}{c}{$0.04019$} & 
\multicolumn{1}{c}{$0.04021$} & \multicolumn{1}{c}{$0.04020$} & 
\multicolumn{1}{c}{$0.04021$} \\ 
&  & \multicolumn{1}{c}{} & \multicolumn{1}{c}{} & \multicolumn{1}{c}{} & 
\multicolumn{1}{c}{} & \multicolumn{1}{c}{} \\ 
0.75 & 0.4 & \multicolumn{1}{c}{$0.62537$} & \multicolumn{1}{c}{$0.62538$} & 
\multicolumn{1}{c}{$0.62543$} & \multicolumn{1}{c}{$0.62541$} & 
\multicolumn{1}{c}{$0.62544$} \\ 
& 0.6 & \multicolumn{1}{c}{$0.48714$} & \multicolumn{1}{c}{$0.48715$} & 
\multicolumn{1}{c}{$0.48723$} & \multicolumn{1}{c}{$0.48719$} & 
\multicolumn{1}{c}{$0.48721$} \\ 
& 0.8 & \multicolumn{1}{c}{$0.37385$} & \multicolumn{1}{c}{$0.37385$} & 
\multicolumn{1}{c}{$0.37394$} & \multicolumn{1}{c}{$0.37390$} & 
\multicolumn{1}{c}{$0.37392$} \\ 
& 1.0 & \multicolumn{1}{c}{$0.28741$} & \multicolumn{1}{c}{$0.28741$} & 
\multicolumn{1}{c}{$0.28750$} & \multicolumn{1}{c}{$0.28746$} & 
\multicolumn{1}{c}{$0.28747$} \\ 
& 3.0 & \multicolumn{1}{c}{$0.02976$} & \multicolumn{1}{c}{$0.02976$} & 
\multicolumn{1}{c}{$0.02978$} & \multicolumn{1}{c}{$0.02977$} & 
\multicolumn{1}{c}{$0.02977$}%
\end{tabular}%
$ \\ \hline
\end{tabular}%
\end{equation*}%
\begin{equation*}
\begin{tabular}{|c|}
\hline
$\text{Table 2 c: Comparison of the numerical solutions of Problem 1 for }%
\lambda =0.01$ \\ \hline
and $N=40,$ $\Delta t=0.0001$ at different times with the exact solutions \\ 
\hline
$%
\begin{tabular}{lllllll}
$x$ & $t$ & $%
\begin{array}{c}
\text{Present} \\ 
\text{ }%
\end{array}%
$ & $%
\begin{array}{c}
\text{Ref.\cite{bs}} \\ 
\text{ }%
\end{array}%
$ & $%
\begin{array}{c}
\text{Ref.\cite{ofem2}} \\ 
\text{ }%
\end{array}%
$ & $%
\begin{array}{c}
\text{Ref.\cite{bbs}} \\ 
\text{ }%
\end{array}%
$ & $%
\begin{array}{c}
\text{Exact} \\ 
\text{ }%
\end{array}%
$ \\ \hline
0.25 & 0.4 & \multicolumn{1}{c}{$0.34192$} & \multicolumn{1}{c}{$0.34192$} & 
\multicolumn{1}{c}{$0.34192$} & \multicolumn{1}{c}{$0.34192$} & 
\multicolumn{1}{c}{$0.34191$} \\ 
& 0.6 & \multicolumn{1}{c}{$0.26897$} & \multicolumn{1}{c}{$0.26897$} & 
\multicolumn{1}{c}{$0.22894$} & \multicolumn{1}{c}{$0.26897$} & 
\multicolumn{1}{c}{$0.22896$} \\ 
& 0.8 & \multicolumn{1}{c}{$0.22148$} & \multicolumn{1}{c}{$0.22148$} & 
\multicolumn{1}{c}{$0.22144$} & \multicolumn{1}{c}{$0.22148$} & 
\multicolumn{1}{c}{$0.22148$} \\ 
& 1.0 & \multicolumn{1}{c}{$0.18819$} & \multicolumn{1}{c}{$0.18819$} & 
\multicolumn{1}{c}{$0.18815$} & \multicolumn{1}{c}{$0.18819$} & 
\multicolumn{1}{c}{$0.18819$} \\ 
& 3.0 & \multicolumn{1}{c}{$0.07511$} & \multicolumn{1}{c}{$0.07511$} & 
\multicolumn{1}{c}{$0.07509$} & \multicolumn{1}{c}{$0.07511$} & 
\multicolumn{1}{c}{$0.07511$} \\ 
&  & \multicolumn{1}{c}{} & \multicolumn{1}{c}{} & \multicolumn{1}{c}{} & 
\multicolumn{1}{c}{} & \multicolumn{1}{c}{} \\ 
0.50 & 0.4 & \multicolumn{1}{c}{$0.66071$} & \multicolumn{1}{c}{$0.66071$} & 
\multicolumn{1}{c}{$0.66075$} & \multicolumn{1}{c}{$0.66071$} & 
\multicolumn{1}{c}{$0.66071$} \\ 
& 0.6 & \multicolumn{1}{c}{$0.52942$} & \multicolumn{1}{c}{$0.52942$} & 
\multicolumn{1}{c}{$0.52946$} & \multicolumn{1}{c}{$0.52942$} & 
\multicolumn{1}{c}{$0.52942$} \\ 
& 0.8 & \multicolumn{1}{c}{$0.43914$} & \multicolumn{1}{c}{$0.43914$} & 
\multicolumn{1}{c}{$0.43917$} & \multicolumn{1}{c}{$0.43914$} & 
\multicolumn{1}{c}{$0.43914$} \\ 
& 1.0 & \multicolumn{1}{c}{$0.37442$} & \multicolumn{1}{c}{$0.37442$} & 
\multicolumn{1}{c}{$0.37444$} & \multicolumn{1}{c}{$0.37442$} & 
\multicolumn{1}{c}{$0.37442$} \\ 
& 3.0 & \multicolumn{1}{c}{$0.15018$} & \multicolumn{1}{c}{$0.15018$} & 
\multicolumn{1}{c}{$0.15016$} & \multicolumn{1}{c}{$0.15018$} & 
\multicolumn{1}{c}{$0.15018$} \\ 
&  & \multicolumn{1}{c}{} & \multicolumn{1}{c}{} & \multicolumn{1}{c}{} & 
\multicolumn{1}{c}{} & \multicolumn{1}{c}{} \\ 
0.75 & 0.4 & \multicolumn{1}{c}{$0.91027$} & \multicolumn{1}{c}{$0.91027$} & 
\multicolumn{1}{c}{$0.91023$} & \multicolumn{1}{c}{$0.91027$} & 
\multicolumn{1}{c}{$0.91026$} \\ 
& 0.6 & \multicolumn{1}{c}{$0.76725$} & \multicolumn{1}{c}{$0.76725$} & 
\multicolumn{1}{c}{$0.76728$} & \multicolumn{1}{c}{$0.76724$} & 
\multicolumn{1}{c}{$0.76724$} \\ 
& 0.8 & \multicolumn{1}{c}{$0.64740$} & \multicolumn{1}{c}{$0.64740$} & 
\multicolumn{1}{c}{$0.64744$} & \multicolumn{1}{c}{$0.64740$} & 
\multicolumn{1}{c}{$0.64740$} \\ 
& 1.0 & \multicolumn{1}{c}{$0.55605$} & \multicolumn{1}{c}{$0.55605$} & 
\multicolumn{1}{c}{$0.55609$} & \multicolumn{1}{c}{$0.55605$} & 
\multicolumn{1}{c}{$0.55605$} \\ 
& 3.0 & \multicolumn{1}{c}{$0.22483$} & \multicolumn{1}{c}{$0.22483$} & 
\multicolumn{1}{c}{$0.22481$} & \multicolumn{1}{c}{$0.22481$} & 
\multicolumn{1}{c}{$0.22481$}%
\end{tabular}%
$ \\ \hline
\end{tabular}%
\end{equation*}

\begin{equation*}
\begin{tabular}{|c|}
\hline
$\text{Table 3 a: Problem 1 for }\lambda =1,$ $t=0.1,$ $\Delta t=0.0001$at
different size with the exact solutions \\ \hline
$%
\begin{tabular}{llllllll}
$%
\begin{array}{c}
\\ 
\begin{array}{c}
\text{Present} \\ 
p=1%
\end{array}%
\end{array}%
$ & $%
\begin{array}{c}
h=0.05 \\ 
\begin{array}{c}
\text{Ref.\cite{bs}} \\ 
\text{ }%
\end{array}%
\end{array}%
$ & $%
\begin{array}{c}
\\ 
\begin{array}{c}
\text{Ref.\cite{bbs}} \\ 
\text{ }%
\end{array}%
\end{array}%
$ &  & $%
\begin{array}{c}
\\ 
\begin{array}{c}
\text{Present} \\ 
p=1%
\end{array}%
\end{array}%
$ & $%
\begin{array}{c}
h=0.025 \\ 
\begin{array}{c}
\text{Ref.\cite{bs}} \\ 
\text{ }%
\end{array}%
\end{array}%
$ & $%
\begin{array}{c}
\\ 
\begin{array}{c}
\text{Ref.\cite{bbs}} \\ 
\text{ }%
\end{array}%
\end{array}%
$ & $%
\begin{array}{c}
\\ 
\begin{array}{c}
\text{Exact} \\ 
\text{ }%
\end{array}%
\end{array}%
$ \\ \hline
\multicolumn{1}{c}{$0.10936$} & \multicolumn{1}{c}{$0.10937$} & 
\multicolumn{1}{c}{$0.10953$} & \multicolumn{1}{c}{} & \multicolumn{1}{c}{$%
0.10949$} & \multicolumn{1}{c}{$0.10949$} & \multicolumn{1}{c}{$0.10954$} & 
\multicolumn{1}{c}{$0.10954$} \\ 
\multicolumn{1}{c}{$0.20943$} & \multicolumn{1}{c}{$0.20945$} & 
\multicolumn{1}{c}{$0.20978$} & \multicolumn{1}{c}{} & \multicolumn{1}{c}{$%
0.20970$} & \multicolumn{1}{c}{$0.20949$} & \multicolumn{1}{c}{$0.20979$} & 
\multicolumn{1}{c}{$0.20979$} \\ 
\multicolumn{1}{c}{$0.29136$} & \multicolumn{1}{c}{$0.29138$} & 
\multicolumn{1}{c}{$0.29188$} & \multicolumn{1}{c}{} & \multicolumn{1}{c}{$%
0.29176$} & \multicolumn{1}{c}{$0.29175$} & \multicolumn{1}{c}{$0.29189$} & 
\multicolumn{1}{c}{$0.29190$} \\ 
\multicolumn{1}{c}{$0.34723$} & \multicolumn{1}{c}{$0.34726$} & 
\multicolumn{1}{c}{$0.34791$} & \multicolumn{1}{c}{} & \multicolumn{1}{c}{$%
0.34775$} & \multicolumn{1}{c}{$0.34773$} & \multicolumn{1}{c}{$0.34792$} & 
\multicolumn{1}{c}{$0.34792$} \\ 
\multicolumn{1}{c}{$0.37076$} & \multicolumn{1}{c}{$0.37080$} & 
\multicolumn{1}{c}{$0.37156$} & \multicolumn{1}{c}{} & \multicolumn{1}{c}{$%
0.37137$} & \multicolumn{1}{c}{$0.37136$} & \multicolumn{1}{c}{$0.37157$} & 
\multicolumn{1}{c}{$0.37158$} \\ 
\multicolumn{1}{c}{$0.35819$} & \multicolumn{1}{c}{$0.35823$} & 
\multicolumn{1}{c}{$0.35902$} & \multicolumn{1}{c}{} & \multicolumn{1}{c}{$%
0.35883$} & \multicolumn{1}{c}{$0.35881$} & \multicolumn{1}{c}{$0.35904$} & 
\multicolumn{1}{c}{$0.35905$} \\ 
\multicolumn{1}{c}{$0.30911$} & \multicolumn{1}{c}{$0.30914$} & 
\multicolumn{1}{c}{$0.30988$} & \multicolumn{1}{c}{} & \multicolumn{1}{c}{$%
0.30971$} & \multicolumn{1}{c}{$0.30969$} & \multicolumn{1}{c}{$0.30990$} & 
\multicolumn{1}{c}{$0.30991$} \\ 
\multicolumn{1}{c}{$0.22719$} & \multicolumn{1}{c}{$0.22722$} & 
\multicolumn{1}{c}{$0.22779$} & \multicolumn{1}{c}{} & \multicolumn{1}{c}{$%
0.22766$} & \multicolumn{1}{c}{$0.22765$} & \multicolumn{1}{c}{$0.22781$} & 
\multicolumn{1}{c}{$0.22782$} \\ 
\multicolumn{1}{c}{$0.12034$} & \multicolumn{1}{c}{$0.12036$} & 
\multicolumn{1}{c}{$0.12066$} & \multicolumn{1}{c}{} & \multicolumn{1}{c}{$%
0.12060$} & \multicolumn{1}{c}{$0.12060$} & \multicolumn{1}{c}{$0.12068$} & 
\multicolumn{1}{c}{$0.12069$}%
\end{tabular}%
$ \\ \hline
\end{tabular}%
\end{equation*}

\begin{equation*}
\begin{tabular}{|c|}
\hline
$\text{Table 3 b: Problem 1 for }\lambda =1,$ $t=0.1,$ $\Delta t=0.0001$at
different size with the exact solutions \\ \hline
$%
\begin{tabular}{llllllll}
$%
\begin{array}{c}
\\ 
\begin{array}{c}
\text{Present} \\ 
p=1%
\end{array}%
\end{array}%
$ & $%
\begin{array}{c}
h=0.0125 \\ 
\begin{array}{c}
\text{Ref.\cite{bs}} \\ 
\text{ }%
\end{array}%
\end{array}%
$ & $%
\begin{array}{c}
\\ 
\begin{array}{c}
\text{Ref.\cite{bbs}} \\ 
\text{ }%
\end{array}%
\end{array}%
$ &  & $%
\begin{array}{c}
\\ 
\begin{array}{c}
\text{Present} \\ 
p=1%
\end{array}%
\end{array}%
$ & $%
\begin{array}{c}
h=0.00625 \\ 
\begin{array}{c}
\text{Ref.\cite{bs}} \\ 
\text{ }%
\end{array}%
\end{array}%
$ & $%
\begin{array}{c}
\\ 
\begin{array}{c}
\text{Ref.\cite{bbs}} \\ 
\text{ }%
\end{array}%
\end{array}%
$ & $%
\begin{array}{c}
\\ 
\begin{array}{c}
\text{Exact} \\ 
\text{ }%
\end{array}%
\end{array}%
$ \\ \hline
\multicolumn{1}{c}{$0.10953$} & \multicolumn{1}{c}{$0.10952$} & 
\multicolumn{1}{c}{$0.10954$} & \multicolumn{1}{c}{} & \multicolumn{1}{c}{$%
0.10954$} & \multicolumn{1}{c}{$0.10953$} & \multicolumn{1}{c}{$0.10954$} & 
\multicolumn{1}{c}{$0.10954$} \\ 
\multicolumn{1}{c}{$0.20977$} & \multicolumn{1}{c}{$0.20975$} & 
\multicolumn{1}{c}{$0.20979$} & \multicolumn{1}{c}{} & \multicolumn{1}{c}{$%
0.20979$} & \multicolumn{1}{c}{$0.20977$} & \multicolumn{1}{c}{$0.20979$} & 
\multicolumn{1}{c}{$0.20979$} \\ 
\multicolumn{1}{c}{$0.29186$} & \multicolumn{1}{c}{$0.29184$} & 
\multicolumn{1}{c}{$0.29189$} & \multicolumn{1}{c}{} & \multicolumn{1}{c}{$%
0.29189$} & \multicolumn{1}{c}{$0.29186$} & \multicolumn{1}{c}{$0.29190$} & 
\multicolumn{1}{c}{$0.29190$} \\ 
\multicolumn{1}{c}{$0.34788$} & \multicolumn{1}{c}{$0.34788$} & 
\multicolumn{1}{c}{$0.34792$} & \multicolumn{1}{c}{} & \multicolumn{1}{c}{$%
0.34792$} & \multicolumn{1}{c}{$0.34788$} & \multicolumn{1}{c}{$0.34792$} & 
\multicolumn{1}{c}{$0.34792$} \\ 
\multicolumn{1}{c}{$0.37153$} & \multicolumn{1}{c}{$0.37153$} & 
\multicolumn{1}{c}{$0.37158$} & \multicolumn{1}{c}{} & \multicolumn{1}{c}{$%
0.37156$} & \multicolumn{1}{c}{$0.37153$} & \multicolumn{1}{c}{$0.37158$} & 
\multicolumn{1}{c}{$0.37158$} \\ 
\multicolumn{1}{c}{$0.35899$} & \multicolumn{1}{c}{$0.35896$} & 
\multicolumn{1}{c}{$0.35904$} & \multicolumn{1}{c}{} & \multicolumn{1}{c}{$%
0.35903$} & \multicolumn{1}{c}{$0.35900$} & \multicolumn{1}{c}{$0.35904$} & 
\multicolumn{1}{c}{$0.35905$} \\ 
\multicolumn{1}{c}{$0.30986$} & \multicolumn{1}{c}{$0.30983$} & 
\multicolumn{1}{c}{$0.30990$} & \multicolumn{1}{c}{} & \multicolumn{1}{c}{$%
0.30989$} & \multicolumn{1}{c}{$0.30986$} & \multicolumn{1}{c}{$0.30990$} & 
\multicolumn{1}{c}{$0.30991$} \\ 
\multicolumn{1}{c}{$0.22778$} & \multicolumn{1}{c}{$0.22776$} & 
\multicolumn{1}{c}{$0.22782$} & \multicolumn{1}{c}{} & \multicolumn{1}{c}{$%
0.22781$} & \multicolumn{1}{c}{$0.22778$} & \multicolumn{1}{c}{$0.22782$} & 
\multicolumn{1}{c}{$0.22782$} \\ 
\multicolumn{1}{c}{$0.12067$} & \multicolumn{1}{c}{$0.12065$} & 
\multicolumn{1}{c}{$0.12069$} & \multicolumn{1}{c}{} & \multicolumn{1}{c}{$%
0.12068$} & \multicolumn{1}{c}{$0.12067$} & \multicolumn{1}{c}{$0.12069$} & 
\multicolumn{1}{c}{$0.12069$}%
\end{tabular}%
$ \\ \hline
\end{tabular}%
\end{equation*}

\textbf{(b)}\ As the second example, we consider particular solution of
Burgers' equation with initial condition

\begin{equation*}
U(x,1)=\exp (\frac{1}{8\lambda }),0\leq x\leq 1,
\end{equation*}%
and boundary conditions $U(0,t)=0$ and $U(1,t)=0$.

This problem has the following analytical solution%
\begin{equation}
U(x,t)=\dfrac{\frac{x}{t}}{1+\sqrt{\frac{t}{t_{0}}}\exp (\frac{x^{2}}{%
4\lambda t})},\text{\quad }t\geq 1,\text{\quad }0\leq x\leq 1,  \label{18}
\end{equation}%
This solution represents the propagation of the shock and the selection of
the smaller\ $\lambda $ result in steep shock solution. So the success of
the numerical method depends on dealing with the steep shock efficiently.

The propogation of the shock is studied with parameters $\lambda
=0.005,0.0005.$ Numerical solutions obtained by exponential collocation
method can be favorably compared with results reported in the papers \cite%
{bs,bbs} at some times in the same Table 4-5. Figs. 10 and 11 show
propagation of shock. $\ $As time advances, the initial steep shock becomes
smoother when the the larger viscosity is used but for the small viscosity
it is steeper. These observations are in complete agreement with those
reported in the papers \cite{nu3}.

\begin{equation*}
\rotatebox{90}{ \begin{tabular}{|l|} \hline $\text{Table 5: Comparison of
results at different times.}\lambda =0.005${\small \ }$[a,b]=[0.1]${\small \
with }$h=0.005$ \\ \hline {\small and }$\Delta t=0.01$ \\ \hline
$\begin{array}{cccccccccc} {\small x} & {\small t=1.7} & {\small t=1.7} &
{\small t=1.7} & {\small t=2.4} & {\small t=2.4} & {\small t=2.4} & {\small
t=3.1} & {\small t=3.1} & {\small t=3.1} \\ & \text{{\small Present}} &
\begin{array}{c} \text{{\small Ref. \cite{bbs}}} \\ \text{{\small
(CBGM)}}\end{array} & \text{{\small Exact}} & \text{{\small Present}} &
\begin{array}{c} \text{{\small Ref. \cite{bbs}}} \\ \text{{\small
(CBGM)}}\end{array} & \text{{\small Exact}} & \text{{\small Present}} &
\begin{array}{c} \text{{\small Ref. \cite{bbs}}} \\ \text{{\small
(CBGM)}}\end{array} & \text{{\small Exact}} \\ {\small 0.1} & {\small
0.058823} & {\small 0.058823} & {\small 0.058823} & {\small 0.041666} &
{\small 0.041666} & {\small 0.041666} & {\small 0.032258} & {\small
0.032258} & {\small 0.032258} \\ {\small 0.2} & {\small 0.117645} & {\small
0.117644} & {\small 0.117645} & {\small 0.083332} & {\small 0.083332} &
{\small 0.083332} & {\small 0.064515} & {\small 0.064515} & {\small
0.064515} \\ {\small 0.3} & {\small 0.176458} & {\small 0.176458} & {\small
0.176458} & {\small 0.124995} & {\small 0.124995} & {\small 0.124995} &
{\small 0.096771} & {\small 0.096771} & {\small 0.096771} \\ {\small 0.4} &
{\small 0.235169} & {\small 0.235170} & {\small 0.235168} & {\small
0.166640} & {\small 0.166639} & {\small 0.166640} & {\small 0.129021} &
{\small 0.129021} & {\small 0.129021} \\ {\small 0.5} & {\small 0.291919} &
{\small 0.291909} & {\small 0.291904} & {\small 0.208115} & {\small
0.208115} & {\small 0.208114} & {\small 0.161231} & {\small 0.161231} &
{\small 0.161231} \\ {\small 0.6} & {\small 0.295900} & {\small 0.294962} &
{\small 0.294910} & {\small 0.247431} & {\small 0.247402} & {\small
0.247417} & {\small 0.193130} & {\small 0.193127} & {\small 0.193127} \\
{\small 0.7} & {\small 0.041931} & {\small 0.042942} & {\small 0.041929} &
{\small 0.252172} & {\small 0.251663} & {\small 0.252172} & {\small
0.221880} & {\small 0.221836} & {\small 0.221867} \\ {\small 0.8} & {\small
0.000639} & {\small 0.000669} & {\small 0.000646} & {\small 0.073049} &
{\small 0.073810} & {\small 0.073025} & {\small 0.215125} & {\small
0.214751} & {\small 0.215135} \\ {\small 0.9} & {\small 0.000005} & {\small
0.000005} & {\small 0.000005} & {\small 0.003008} & {\small 0.003114} &
{\small 0.003023} & {\small 0.070730} & {\small 0.071252} & {\small
0.070874} \\ & & & & & & & & & \\ {\small L}_{2}{\small \times 10}^{3} &
{\small 0.00459} & {\small 0.35126} & & {\small 0.00457} & {\small 0.24448}
& & {\small 0.02295} & {\small 0.235} & \\ {\small L}_{\infty }{\small
\times 10}^{3} & {\small 0.06494} & {\small 1.20726} & & {\small 0.04639} &
{\small 0.80176} & & {\small 0.42187} & {\small 4.79061} & \end{array}$ \\
\hline \end{tabular}}
\end{equation*}

\begin{equation*}
\vspace{0.2cm}%
\begin{tabular}{|l|}
\hline
$\text{Table 4: Comparison of results at different times.}\lambda =0.0005$
with $h=0.005$ and $\Delta t=0.01$ \\ \hline
$%
\begin{array}{cccccccccc}
x & t=1.7 & t=1.7 & t=1.7 & t=2.5 & t=2.5 & t=2.5 & t=3.25 & t=3.25 & t=3.25
\\ 
& \text{Present} & \text{Ref. \cite{bs}} & \text{Exact} & \text{Present} & 
\text{Ref. \cite{bs}} & \text{Exact} & \text{Present} & \text{Ref. \cite{bs}}
& \text{Exact} \\ 
{\small 0.1} & {\small 0.05882} & {\small 0.05883} & {\small 0.05882} & 
{\small 0.04000} & {\small 0.04000} & {\small 0.04000} & {\small 0.03077} & 
{\small 0.03077} & {\small 0.03077} \\ 
{\small 0.2} & {\small 0.11765} & {\small 0.11765} & {\small 0.11765} & 
{\small 0.08000} & {\small 0.08000} & {\small 0.08000} & {\small 0.06154} & 
{\small 0.06154} & {\small 0.06154} \\ 
{\small 0.3} & {\small 0.17647} & {\small 0.17648} & {\small 0.17647} & 
{\small 0.12000} & {\small 0.12001} & {\small 0.12000} & {\small 0.09231} & 
{\small 0.09231} & {\small 0.09231} \\ 
{\small 0.4} & {\small 0.23529} & {\small 0.23531} & {\small 0.23529} & 
{\small 0.16000} & {\small 0.16001} & {\small 0.16000} & {\small 0.12308} & 
{\small 0.12308} & {\small 0.12308} \\ 
{\small 0.5} & {\small 0.29412} & {\small 0.29414} & {\small 0.29412} & 
{\small 0.20000} & {\small 0.20001} & {\small 0.20000} & {\small 0.15385} & 
{\small 0.15385} & {\small 0.15385} \\ 
{\small 0.6} & {\small 0.35294} & {\small 0.35296} & {\small 0.35294} & 
{\small 0.24000} & {\small 0.24001} & {\small 0.24000} & {\small 0.18462} & 
{\small 0.18462} & {\small 0.18462} \\ 
{\small 0.7} & {\small 0.00000} & {\small 0.00000} & {\small 0.00000} & 
{\small 0.28000} & {\small 0.28001} & {\small 0.28000} & {\small 0.21538} & 
{\small 0.21539} & {\small 0.21538} \\ 
{\small 0.8} & {\small 0.00000} & {\small 0.00000} & {\small 0.00000} & 
{\small 0.00828} & {\small 0.00811} & {\small 0.00977} & {\small 0.24615} & 
{\small 0.24616} & {\small 0.24615} \\ 
{\small 0.9} & {\small 0.00000} & {\small 0.00000} & {\small 0.00000} & 
{\small 0.00000} & {\small 0.00000} & {\small 0.00000} & {\small 0.12394} & 
{\small 0.12358} & {\small 0.12435}%
\end{array}%
$ \\ \hline
\end{tabular}%
\end{equation*}

\begin{equation*}
\begin{array}{cc}
\begin{tabular}{l}
\FRAME{itbpF}{2.8729in}{2.4275in}{0in}{}{}{fig10.bmp}{\special{language
"Scientific Word";type "GRAPHIC";display "USEDEF";valid_file "F";width
2.8729in;height 2.4275in;depth 0in;original-width 14.2288in;original-height
8.0004in;cropleft "0";croptop "1";cropright "1";cropbottom "0";filename
'fig/Fig10.bmp';file-properties "XNPEU";}} \\ 
Fig. 10: Shock propagation, $\lambda =0.005$%
\end{tabular}
& 
\begin{tabular}{l}
\FRAME{itbpF}{3.0511in}{2.4267in}{0in}{}{}{fig11.bmp}{\special{language
"Scientific Word";type "GRAPHIC";display "USEDEF";valid_file "F";width
3.0511in;height 2.4267in;depth 0in;original-width 14.2288in;original-height
8.0004in;cropleft "0";croptop "1";cropright "1";cropbottom "0";filename
'fig/Fig11.bmp';file-properties "XNPEU";}} \\ 
Fig. 11: Shock propagation, $\lambda =0.0005$%
\end{tabular}%
\end{array}%
\end{equation*}

\textbf{(c)\ }Travelling wave solution of the Burgers' equation has the
form: 
\begin{equation}
U(x,t)=\dfrac{\alpha +\mu +(\mu -\alpha )\exp \eta }{1+\exp \eta },\text{ }%
0\leq x\leq 1,\text{ }t\geq 0,  \label{20}
\end{equation}%
where 
\begin{equation*}
\eta =\dfrac{\alpha (x-\mu t-\gamma )}{\lambda },
\end{equation*}%
and $\alpha ,~\mu $ and $\gamma $ are arbitrary constants. The boundary
conditions is 
\begin{eqnarray*}
U(0,t) &=&1,\text{ }U(1,t)=0.2\text{ or} \\
U_{x}(0,t) &=&0,\text{ }U_{x}(1,t)=0,\text{ for }t\geq 0
\end{eqnarray*}%
and initial condition is obtained from the analytical solution (\ref{20})
when $t=0$. Analytical solution takes values between 1 and 0.2 and the
propagation of the wave front through the right will be observed with
varying $\lambda .$ \ The smaller $\lambda $ we take for the Burgers'
equation , the steeper the wave front propagates \ The robustness of the
algorithm will be shown by monitoring the motion of the wave front with
smaller $\lambda $. The algorithm has run for the values $\alpha =0.4,$ $\mu
=0.6,$ $\gamma =0.125$ and $\lambda =0.01,h=1/36,\Delta t=0.001$. A
comparison of the results obtained from the present method with those by \ I
Dag and his coworkers is shown in table 5. $L_{\infty }$-norm of the methods
given data in the table 6 has been found as $0.004,0.004,0.005,0.006,0.005$
respectively. Numerical solutions obtained with present method are in good
agreement with results obtained by the methods documented in the Table 6.
Visual motion of the wave front is depicted in the Figs. 12-13 for the $%
\lambda =0.01,0.05.$ The numerical results demonstrate the formation of the
steep front and very steeper front. \ Error graphs of the numerical
solutions are also shown in the Figs. 14-15, from figures the maximum error
occurs in the midle of the solution domain. Solutions from time $t=0$ to $%
t=1.2$ at some times are visualised in 3D graph to see the propogation of
the sharp behaviours in Figs 16-17.

\begin{equation*}
\begin{tabular}{|l|}
\hline
Table 6: Comparison of results at time \\ \hline
$t=0.5,$ $h=1/36,$ $\Delta t=0.01,$ $\lambda =0.01$ \\ \hline
$%
\begin{array}{lccccc}
x & 
\begin{array}{c}
\text{Present} \\ 
p=1%
\end{array}
& 
\begin{array}{c}
\text{Ref. \cite{bs}} \\ 
\Delta t=0.025%
\end{array}
& 
\begin{array}{c}
\text{Ref. \cite{bbs}} \\ 
\text{(QBGM)}%
\end{array}
& 
\begin{array}{c}
\text{Ref. \cite{bbs}} \\ 
\text{(CBGM)}%
\end{array}
& 
\begin{array}{c}
\text{Exact} \\ 
\end{array}
\\ 
{\small 0.000} & {\small 1.} & {\small 1.} & {\small 1.} & {\small 1.} & 
{\small 1.} \\ 
{\small 0.056} & {\small 1.} & {\small 1.} & {\small 1.} & {\small 1.} & 
{\small 1.} \\ 
{\small 0.111} & {\small 1.} & {\small 1.} & {\small 1.} & {\small 1.} & 
{\small 1.} \\ 
{\small 0.167} & {\small 1.} & {\small 1.} & {\small 1.} & {\small 1.} & 
{\small 1.} \\ 
{\small 0.222} & {\small 1.} & {\small 1.} & {\small 1.} & {\small 1.} & 
{\small 1.} \\ 
{\small 0.278} & {\small 0.998} & {\small 0.999} & {\small 0.998} & {\small %
0.998} & {\small 0.998} \\ 
{\small 0.333} & {\small 0.980} & {\small 0.986} & {\small 0.980} & {\small %
0.980} & {\small 0.980} \\ 
{\small 0.389} & {\small 0.847} & {\small 0.850} & {\small 0.841} & {\small %
0.842} & {\small 0.847} \\ 
{\small 0.444} & {\small 0.452} & {\small 0.448} & {\small 0.458} & {\small %
0.457} & {\small 0.452} \\ 
{\small 0.500} & {\small 0.238} & {\small 0.236} & {\small 0.240} & {\small %
0.241} & {\small 0.238} \\ 
{\small 0.556} & {\small 0.204} & {\small 0.204} & {\small 0.205} & {\small %
0.205} & {\small 0.204} \\ 
{\small 0.611} & {\small 0.2} & {\small 0.2} & {\small 0.2} & {\small 0.2} & 
{\small 0.2} \\ 
{\small 0.667} & {\small 0.2} & {\small 0.2} & {\small 0.2} & {\small 0.2} & 
{\small 0.2} \\ 
{\small 0.722} & {\small 0.2} & {\small 0.2} & {\small 0.2} & {\small 0.2} & 
{\small 0.2} \\ 
{\small 0.778} & {\small 0.2} & {\small 0.2} & {\small 0.2} & {\small 0.2} & 
{\small 0.2} \\ 
{\small 0.833} & {\small 0.2} & {\small 0.2} & {\small 0.2} & {\small 0.2} & 
{\small 0.2} \\ 
{\small 0.889} & {\small 0.2} & {\small 0.2} & {\small 0.2} & {\small 0.2} & 
{\small 0.2} \\ 
{\small 0.944} & {\small 0.2} & {\small 0.2} & {\small 0.2} & {\small 0.2} & 
{\small 0.2} \\ 
{\small 1.000} & {\small 0.2} & {\small 0.2} & {\small 0.2} & {\small 0.2} & 
{\small 0.2}%
\end{array}%
$ \\ \hline
\end{tabular}%
\end{equation*}%
\begin{equation*}
\begin{array}{c}
\FRAME{itbpF}{2.8738in}{2.4275in}{0in}{}{}{fig12.bmp}{\special{language
"Scientific Word";type "GRAPHIC";display "USEDEF";valid_file "F";width
2.8738in;height 2.4275in;depth 0in;original-width 14.2288in;original-height
8.0004in;cropleft "0";croptop "1";cropright "1";cropbottom "0";filename
'fig/Fig12.bmp';file-properties "XNPEU";}} \\ 
\multicolumn{1}{l}{\text{Figure 12: Solutions at different times for }%
\lambda =0.01} \\ 
\multicolumn{1}{l}{h=1/36,\text{ }\Delta t=0.001,\text{ }x\in \lbrack 0,1].}%
\end{array}%
\begin{tabular}{l}
\FRAME{itbpF}{2.8738in}{2.4275in}{0in}{}{}{fig13.bmp}{\special{language
"Scientific Word";type "GRAPHIC";display "USEDEF";valid_file "F";width
2.8738in;height 2.4275in;depth 0in;original-width 14.2288in;original-height
8.0004in;cropleft "0";croptop "1";cropright "1";cropbottom "0";filename
'fig/Fig13.bmp';file-properties "XNPEU";}} \\ 
Figure 13: $L_{2}$ error norm for $\lambda =0.005$ \\ 
$h=1/36,\text{ }p=1,$ $\Delta t=0.001,\text{ }x\in \lbrack 0,1].$%
\end{tabular}%
\end{equation*}%
\begin{equation*}
\begin{array}{cc}
\begin{tabular}{l}
\FRAME{itbpF}{2.9594in}{2.3687in}{0in}{}{}{fig14.bmp}{\special{language
"Scientific Word";type "GRAPHIC";display "USEDEF";valid_file "F";width
2.9594in;height 2.3687in;depth 0in;original-width 4.1874in;original-height
3.3434in;cropleft "0";croptop "1";cropright "1";cropbottom "0";filename
'fig/Fig14.bmp';file-properties "XNPEU";}} \\ 
$\text{Fig. }$14: $L_{2}$ error norm for $\lambda =0.01$ and $h=1/36$%
\end{tabular}
& 
\begin{tabular}{l}
\FRAME{itbpF}{2.9014in}{2.4933in}{0in}{}{}{fig15.bmp}{\special{language
"Scientific Word";type "GRAPHIC";display "USEDEF";valid_file "F";width
2.9014in;height 2.4933in;depth 0in;original-width 4.1044in;original-height
3.5206in;cropleft "0";croptop "1";cropright "1";cropbottom "0";filename
'fig/Fig15.bmp';file-properties "XNPEU";}} \\ 
$\text{Fig. 15 }L_{2}$ error norm $\lambda =0.005$ and $h=1/36$%
\end{tabular}%
\end{array}%
\end{equation*}%
\begin{equation*}
\begin{array}{cc}
\begin{tabular}{l}
\FRAME{itbpF}{3.0441in}{2.879in}{0in}{}{}{fig16.bmp}{\special{language
"Scientific Word";type "GRAPHIC";maintain-aspect-ratio TRUE;display
"USEDEF";valid_file "F";width 3.0441in;height 2.879in;depth
0in;original-width 5.0272in;original-height 4.7539in;cropleft "0";croptop
"1";cropright "1";cropbottom "0";filename 'fig/Fig16.bmp';file-properties
"XNPEU";}} \\ 
Fig. 16: Shock propagation, $\lambda =0.01$%
\end{tabular}
& 
\begin{tabular}{l}
\FRAME{itbpF}{3.1514in}{2.8426in}{0.1565in}{}{}{fig17.bmp}{\special{language
"Scientific Word";type "GRAPHIC";maintain-aspect-ratio TRUE;display
"USEDEF";valid_file "F";width 3.1514in;height 2.8426in;depth
0.1565in;original-width 5.2062in;original-height 4.6933in;cropleft
"0";croptop "1";cropright "1";cropbottom "0";filename
'fig/Fig17.bmp';file-properties "XNPEU";}} \\ 
Fig. 17: Shock propagation, $\lambda =0.005$%
\end{tabular}%
\end{array}%
\end{equation*}%
\ 

\textbf{(d)} The initial condition and boundary conditions are taken as%
\begin{equation*}
U(x,0)=\lambda \left[ x+\tan (x/2)\right] ,
\end{equation*}

\begin{equation*}
\begin{tabular}{l}
$U(0.5,t)=\dfrac{\lambda }{1+\lambda t}\left[ 0.5+\tan \left( \dfrac{1}{%
4(1+\lambda t)}\right) \right] ,\quad t\geq 0,\vspace{0.2cm}$ \\ 
$U(1.5,t)=\dfrac{\lambda }{1+\lambda t}\left[ 1.5+\tan \left( \dfrac{3}{%
4(1+\lambda t)}\right) \right] ,$\quad $t\geq 0,$%
\end{tabular}%
\end{equation*}%
The exact solution of the problem is \cite{nu2}%
\begin{equation*}
U(x,t)=\frac{\lambda }{1+\lambda t}\left[ x+\tan \left( \frac{x}{2(1+\lambda
t)}\right) \right] .
\end{equation*}

Table 7 gives results of the problem (d) among exponential cubic B-spline
collocation method, cubic B-spline collocation method and exact solutions
for $\lambda =$ $1/1000$ at time $t=2.25$.$\vspace{0.2cm}$The program is run
for a set of parameters $p$, the best results is found for the value $p=1$
seeing in the Table 7.

\hspace{-2.5cm}%
\begin{equation*}
\begin{tabular}{|l|}
\hline
Table 7: Comparison of results at time \\ \hline
\begin{tabular}{l}
$t=2.25,$ $\lambda =1/1000$ with $h=0.00625,$ $\Delta t=0.015$ \\ \hline
$%
\begin{array}{c|cccccc}
x & 
\begin{array}{c}
\text{Present} \\ 
(p=0.01)%
\end{array}
& 
\begin{array}{c}
\text{Present} \\ 
(p=0.1)%
\end{array}
& 
\begin{array}{c}
\text{Present} \\ 
(p=0.5)%
\end{array}
& 
\begin{array}{c}
\text{Present} \\ 
(p=1)%
\end{array}
& 
\begin{array}{c}
\text{cubic B-spline} \\ 
\text{\cite{bs}}%
\end{array}
& \text{Exact} \\ \hline
0.5 & {\small 0.007329977} & {\small 0.007329977} & {\small 0.000753049} & 
{\small 0.000753049} & {\small 0.007329977} & {\small 0.000753049} \\ 
0.6 & {\small 0.008817897} & {\small 0.008817897} & {\small 0.000905597} & 
{\small 0.000905597} & {\small 0.008817897} & {\small 0.000906558} \\ 
0.7 & {\small 0.010323604} & {\small 0.010323604} & {\small 0.001061678} & 
{\small 0.001061678} & {\small 0.010323604} & {\small 0.001061749} \\ 
0.8 & {\small 0.011851020} & {\small 0.011851020} & {\small 0.001218899} & 
{\small 0.001218991} & {\small 0.011851020} & {\small 0.001218992} \\ 
0.9 & {\small 0.013402144} & {\small 0.013402144} & {\small 0.001378707} & 
{\small 0.001378707} & {\small 0.013402144} & {\small 0.001378707} \\ 
1.0 & {\small 0.014978854} & {\small 0.014978854} & {\small 0.001541377} & 
{\small 0.001541377} & {\small 0.014978854} & {\small 0.001541377} \\ 
1.1 & {\small 0.016584692} & {\small 0.016584692} & {\small 0.001707565} & 
{\small 0.001707565} & {\small 0.016584692} & {\small 0.001707565} \\ 
1.2 & {\small 0.016584692} & {\small 0.018226288} & {\small 0.001877935} & 
{\small 0.001877935} & {\small 0.018226288} & {\small 0.001877935} \\ 
1.3 & {\small 0.019914461} & {\small 0.019914461} & {\small 0.002053284} & 
{\small 0.002053119} & {\small 0.019914461} & {\small 0.002053284} \\ 
1.4 & {\small 0.021662973} & {\small 0.021662973} & {\small 0.002231955} & 
{\small 0.002231955} & {\small 0.021662973} & {\small 0.002234577} \\ 
1.5 & {\small 0.023483948} & {\small 0.023483948} & {\small 0.002423004} & 
{\small 0.002423004} & {\small 0.023483948} & {\small 0.002423004}%
\end{array}%
$%
\end{tabular}
\\ \hline
\end{tabular}%
\end{equation*}%
Visual representation of the absolute error over space interval $0.5\leq
x\leq 1.5$ is given at time $t=2.25$ in Figs18-19.%
\begin{equation*}
\begin{tabular}{l}
\FRAME{itbpF}{2.8738in}{2.4275in}{0in}{}{}{fig18.bmp}{\special{language
"Scientific Word";type "GRAPHIC";display "USEDEF";valid_file "F";width
2.8738in;height 2.4275in;depth 0in;original-width 14.2288in;original-height
8.0004in;cropleft "0";croptop "1";cropright "1";cropbottom "0";filename
'fig/Fig18.bmp';file-properties "XNPEU";}} \\ 
Figure 18: Absolute error for $h=1/36,$ \\ 
$p=1,$ $\Delta t=0.001,$ $t=2.25,$ $\lambda =0.1$%
\end{tabular}%
\begin{tabular}{l}
\FRAME{itbpF}{2.8738in}{2.4275in}{0in}{}{}{fig19.bmp}{\special{language
"Scientific Word";type "GRAPHIC";display "USEDEF";valid_file "F";width
2.8738in;height 2.4275in;depth 0in;original-width 14.2288in;original-height
8.0004in;cropleft "0";croptop "1";cropright "1";cropbottom "0";filename
'fig/Fig19.bmp';file-properties "XNPEU";}} \\ 
Figure 19: Absolute error for $h=1/36,$ \\ 
$p=1,$ $\Delta t=0.001,$ $t=2.25,$ $\lambda =0.01$%
\end{tabular}%
\end{equation*}

In Table 8 the accuracy of the present method via $L_{\infty }$ norm is can
be examined. As number of $N$ and $\lambda $ parameter increase, according
to Table 8, the error decrease.%
\begin{equation*}
\begin{tabular}{|l|}
\hline
Table 8: Comparison of results at time $t=2.25$ max error. \\ \hline
\begin{tabular}{l}
Parameters $p=0.1,$ $\Delta t=0.01$ \\ \hline
$%
\begin{array}{ccccc}
\lambda & N=16 & N=32 & N=64 & N=128 \\ 
1 & 1.19\times 10^{-7} & 5.26\times 10^{-8} & 3.94\times 10^{-8} & 
3.73\times 10^{-8} \\ 
1/2 & 1.03\times 10^{-7} & 4.69\times 10^{-8} & 2.15\times 10^{-8} & 
1.48\times 10^{-8} \\ 
1/4 & 7.19\times 10^{-6} & 2.74\times 10^{-6} & 1.53\times 10^{-6} & 
1.22\times 10^{-6} \\ 
1/8 & 6.57\times 10^{-5} & 2.06\times 10^{-5} & 9.36\times 10^{-6} & 
6.52\times 10^{-6} \\ 
1/16 & 1.45\times 10^{-4} & 4.09\times 10^{-5} & 1.51\times 10^{-5} & 
8.65\times 10^{-6} \\ 
1/32 & 1.58\times 10^{-4} & 4.18\times 10^{-5} & 1.31\times 10^{-5} & 
5.98\times 10^{-6} \\ 
1/100 & 1.32\times 10^{-4} & 3.40\times 10^{-5} & 9.2\times 10^{-6} & 
3.04\times 10^{-6}%
\end{array}%
$%
\end{tabular}
\\ \hline
\end{tabular}%
\end{equation*}

\section{Conclusion}

The Exponential cubic B-Spline collocation method for the numerical
solutions of the Burges' equation is presented over the finite elements so
that the continuity of the dependent variable and its first two derivatives
is satisfied \ for the approximate solution throughout the solution range.
The equation has been integrated into a system of the linearized iterative
algebraic equations. The system of the iterative at each time step in which
it has got a three-banded coefficients matrix\ is solved with the Thomas
algorithm. Generally, Comparative results show that results of our finding
is better than the that of the cubic B-spline collocation method and is much
the same with that of the cubic B-spline Galerkin finite element method.
Since cost of the cubic B-spline Galerkin method is higher than the
suggested method, \ that is advantages of the Exponential cubic B-Spline
collocation method over the cubic B-spline Galerkin method. During all runs
of the algorithm, the best result are found for the free parameter $p=1$ for
the Exponential cubic B-Spline functions.

\end{document}